\documentclass[11pt,a4paper,twoside]{amsart}
\usepackage{amsfonts}
\usepackage{amsmath}
\usepackage{amsthm}
\usepackage{amscd}
\usepackage{euscript}

\newcommand{\OO}{{\EuScript O}}

\newcommand{\PP}{{\mathbb P}}

\newcommand{\ZZ}{{\mathbb Z}}

\newcommand{\ch}{\vee}
\newcommand{\Vop}{{ V_0 \otimes \dots \otimes V_p}}
\newcommand{\Woq}{{ W_0 \otimes \dots \otimes W_q}}

\newcommand{\PVop}{{ \PP(V_1) \times \dots \times \PP(V_{p-1})}}
\newcommand{\PWoq}{{ \PP(W_1) \times \dots \times \PP(W_{q-1})}}

\swapnumbers
\newtheorem{thm}{Theorem}[section]

\newtheorem{cor}[thm]{Corollary}
\newtheorem{prop}[thm]{Proposition}
\theoremstyle{definition}
\newtheorem{defin}[thm]{Definition}

\theoremstyle{remark}
\newtheorem{remark}[thm]{Remark}

 \def\rig#1{\smash{ \mathop{\longrightarrow}
    \limits^{#1}}}

\title{\large {The Binet-Cauchy Theorem for the Hyperdeterminant of 
boundary format multidimensional Matrices}}

\begin{document}

\maketitle

\markboth{The Binet-Cauchy Theorem for the Hyperdeterminant of 
Boundary Format Multidimensional Matrices}
{CARLA DIONISI AND  GIORGIO OTTAVIANI}

  \author{\begin{center}
        CARLA DIONISI \\
  \begin{small}
  {\it Dipartimento di Matematica Applicata ``G.~Sansone'', via S.Marta
    3, 50139 Firenze, Italy}\\ 
  {\it e-mail address}: {\tt dionisi@math.unifi.it}\\ 
    \end{small}
  
  GIORGIO OTTAVIANI \footnote{Both the authors    
        were partially supported by Italian MURST and GNSAGA-INDAM}\\
  \begin{small}
   {\it Dipartimento di Matematica ``U.~Dini'', Viale Morgagni 67/a, 
    50134 Firenze, Italy}\\
  {\it e-mail address}: {\tt ottavian@math.unifi.it}
      \end{small}
    \end{center}
  }

  \markboth{The Binet-Cauchy Theorem for the Hyperdeterminant of 
  Boundary Format Multidimenional Matrices}
  {CARLA DIONISI AND  GIORGIO OTTAVIANI}


\section{Introduction} \indent
The Binet-Cauchy Theorem states that if $A$ and $B$ are square
matrices then \qquad $det(A\cdot B)=det(A)\cdot det(B)$. The main result of
this paper is a generalization of this theorem to multidimensional
matrices $A$, $B$ of boundary format (see definition
\ref{def:boundary}), where the hyperdeterminant replaces the
determinant (see the theorem (\ref{binet}) for the precise statement).
The idea of the proof is quite simple, in fact we consider the
hyperdeterminant of $A$ as the determinant of a certain morphism
$\partial_A$ (see definition \ref{defdet}) as in \cite{GKZ}.  Then we compute
$\partial_{A\ast B}$ by means of $\partial_A$ and $\partial_B$ and we
apply the usual Binet-Cauchy Theorem.  The proof is better understood
with the language of vector bundles in the setting of algebraic
geometry, although we do not strictly need them.  The study of
multiplicative properties of hyperdeterminants was left as an open
problem in \cite{GKZ}.
 As a consequence (corollary \ref{nondeg}), we prove that
given two matrices $A$ and $B$ of boundary format then $A \ast B$ is
nondegenerate if and only if $A$ and $B$ are both nondegenerate. We
show by a counterexample (remark \ref{rem:contro}) that the
assumption of boundary format cannot be dropped.

We remind how the definition of hyperdeterminant comes out. In chapter
14 of \cite{GKZ} the hyperdeterminant is defined geometrically by
considering the dual variety, that is by studying tangency conditions.

The nondegeneracy of a multidimensional matrix is algebraically
equivalent to the absence of nontrivial solutions of a suitable system
of equations containing some partial derivatives.  With this approach
the usual determinant of a square $n \times n$ matrix is realized as
the equation of the dual variety to the Segre variety $\PP^{n-1}
\times \PP^{n-1}$.

A second well known approach is to define a square matrix to be
nondegenerate if the associated linear system has only trivial
solutions.  In this paper we choose this second approach as the
definition of nondegeneracy (\ref{def:mat}). The nondegenerate
matrices fill up a codimension one subvariety exactly in the boundary
format case.  In this case the second approach is simpler and it
allows us to compute the degree of the hyperdeterminant and to give an
explicit formula for it directly from this definition of
nondegeneracy. The above results were found in \cite{GKZ} as
consequences of a combinatorial statement (lemma 14.2.7) which needs a
nontrivial proof about the irreducibility over $\ZZ$ of a certain
polynomial (14.3.4, 14.3.5, 14.3.6 of \cite{GKZ}). Following this
approach, theorem 3.3 of \cite{GKZ} comes quickly and the computation
of the degree of the hyperdeterminant is a trivial consequence.

Our definition fits into invariant theory and does not depend on
coordinates. The tools that we use are vector
bundles over the product of projective spaces (as in \cite{AO99} or \cite{D}) and
K\"unneth formula to compute their cohomology.

In the remark \ref{notone} we notice that an analog of the hyperdeterminant can be defined 
also in some cases where the variety of degenerate matrices has big
codimension. This fact seems promising for other applications (see \cite{CO}).

\section{Notations and preliminaries} \indent
Let $V_i$ for $i=0,\ldots ,p$ be a complex vector space of dimension
$k_i+1$. We assume $k_0=\max_{i}k_i$. It is not necessary to assume
$k_0\ge k_1\ge\ldots\ge k_p$ (see remark \ref{perm}). 

We remark that a multidimensional matrix  $A \in V_0 \otimes \ldots
  \otimes V_p$ can be regarded as a map $V_0^ \ch \to V_1 \otimes
  \dots \otimes V_p$, hence taken
  the dual map $V_1^ \ch \otimes \ldots \otimes V_p^\ch \to V_0$ (that
  we call also $A$), we can give the following definition:
\begin{defin}
\label{def:mat}
A multidimensional matrix $A$ is  called
{\em degenerate} if there are $v_i\in V_i^*$, $v_i\neq 0$ for $i=1,\ldots
,p$ such that $A(v_1\otimes\ldots\otimes v_p)=0$.
 \end{defin}
 
 If $p=1$ nondegenerate matrices are exactly the matrices of maximal
 rank.
 
 If $k_0\ge\sum_{i=1}^pk_i$ it is easy to check (\cite{WZ} and also the 
proof
of theorem (\ref{main})) that
 degenerate matrices fill an irreducible variety of codimension
 $k_0-\sum_{i=1}^pk_i+1$. If $k_0<\sum_{i=1}^pk_i$ then all matrices
 in $V_0\otimes\ldots\otimes V_p$ are degenerate.

\begin{defin}
\label{def:boundary}
  If $k_0=\sum_{i=1}^pk_i$ the matrices $A\in V_0\otimes\ldots\otimes
  V_p$ are called {\em of boundary format}.
\end{defin}
\begin{remark}(see for instance \cite{h})
  \label{rem:det}
For a vector space $V$ of dimension $n$ we denote $ det V:=\wedge^n
V$.
We recall that any linear map $\Phi \in Hom(V,W)$ between vector spaces
of the same dimension induces the map $det \Phi \in Hom(det V, det
W)$. If $A$ and $B$ are vector spaces of dimension $a$ and $b$ respectively,
then there are canonical isomorphisms: 
\[det(A \otimes B)\simeq (det A)^{\otimes b} \otimes (det B)^{\otimes a}
  \qquad det (S^k A)\simeq(det A)^{\otimes \binom{a+k-1}{a} }\]
\[\wedge^k A\simeq \wedge^{a-k}A^*\otimes (det A)\]
The above isomorphisms hold also if $A$ and $B$ are replaced by vector
bundles over a variety $X$.
\end{remark}

\section{Hyperdeterminants} \indent
Let $A\in V_0\otimes\ldots\otimes V_p$ be of boundary format and let $m_j=\sum_{i=1}^{j-1}k_i$ with
the convention $m_1=0$.

We remark that the definition of $m_i$ depends on the order we have
chosen among the $k_j$'s (see remark \ref{perm}).

With the above notations the vector spaces $V_0^\ch \otimes
S^{m_1}V_1\otimes\ldots\otimes S^{m_p}V_p$ and $
S^{m_1+1}V_1\otimes\ldots\otimes S^{m_p+1}V_p$ have the same dimension
$N=\frac{(k_0+1)!}{k_1!\ldots k_r!}.$

The following theorem is essentially equivalent to theorem 4.3 and
lemma 4.4  of \cite{gkz1}.
Since we want to make paper self-contained and since our proof of the
irreducibility of the homogeneous polynomial $Det$ does not need any
combinatorial statement as in \cite{GKZ} and \cite{gkz1}, then we
include the proof.
\begin{thm}{\bf (and definition of $\partial_A$).}
\label{main}
Let $k_0=\sum_{i=1}^pk_i$. Then the degenerate matrices fill an
irreducible subvariety of degree $N=\frac{(k_0+1)!}{k_1!\ldots k_r!}$
whose equation is given by the determinant of the natural morphism
\[\partial_A: V_0^\ch \otimes S^{m_1}V_1\otimes\ldots\otimes
S^{m_p}V_p \rig{}
S^{m_1+1}V_1\otimes\ldots\otimes S^{m_p+1}V_p\]
\end{thm}

\begin{proof}
  If $A$ is degenerate then we get $A(v_1\otimes\ldots\otimes v_p)=0$
  for some $v_i\in V_i^*$, $v_i\neq 0$ for $i=1,\ldots ,p$. Then
  $({\partial}_A)^t\left(v_1^{\otimes m_1+1} \otimes\ldots\otimes v_p^{\otimes
      m_p+1} \right)=0$.
  
  Conversely if $A$ is nondegenerate we get a surjective natural map
  of vector bundles over $X=\PP(V_2)\times\ldots\times\PP(V_p)$
\begin{equation*}
  \label{map}
 V_0^\ch \otimes \OO_X\rig{\phi_A}V_1\otimes\OO_X(1,\ldots ,1).  
\end{equation*}

Indeed, by our definition, $\phi_A$ is surjective if and only if $A$ is nondegenerate.

We construct a vector bundle $S$ over
$\PP(V_2)\times\ldots\times\PP(V_p)$ whose dual $S^*$ is the kernel of
$\phi_A$ so that we have the exact sequence
\begin{equation}
  \label{ker}
0\rig{}S^*\rig{}V_0^\ch\otimes \OO\rig{}V_1\otimes\OO(1,\ldots ,1)\rig{}0. 
\end{equation}

After tensoring by $\OO(m_2,\ldots ,m_p)$ and taking cohomology we get
\begin{equation*}
  \label{eq:3}
H^0(S^*(m_2,m_3,\ldots ,m_p))\rig{}V_0^\ch\otimes S^{m_1}V_1\otimes\ldots\otimes   S^{m_p}V_p\rig{{\partial}_A}
S^{m_1+1}V_1\otimes\ldots\otimes S^{m_p+1}V_p
\end{equation*}

and we need to prove
\begin{equation}
  \label{van}
H^0(S^*(m_2,m_3,\ldots ,m_p))=0.
\end{equation}

Let
$d=\dim\left(\PP(V_2)\times\ldots\times\PP(V_p)\right)=\sum_{i=2}^pk_i=m_{p+1}-k_1$.

Since $det(S^*)=\OO(-k_1-1,\ldots ,-k_1-1)$ and $rank~S^*=d$ 
from remark \ref{rem:det} it follows that
\begin{equation}
  \label{isos}
S^*(m_2,m_3,\ldots ,m_p)\simeq\wedge^{d-1}S(-1,-k_1-1+m_3,\ldots,
-k_1-1+m_p) 
\end{equation}

Hence, by taking  the  $(d-1)$-th wedge power of the dual of
the sequence (\ref{ker}), and using  K\"unneth formula 
to calculate the cohomology as in \cite{gkz1}, the result follows.
In order to prove the irreducibility of the subvariety $D$ of degenerate 
matrices it is sufficient to construct the incidence variety
\[Z=\{\left(A,([v_1],\ldots ,[v_p]\right)\in \left(V_0\otimes\ldots
V_p\right)\times\left[\PP(V_1)\times\ldots\times\PP(V_p)\right] |
A(v_1\otimes\ldots \otimes v_p)=0\}\]
$Z$ is a vector bundle over $\PP(V_1)\times\ldots\times\PP(V_p)$,
hence it is irreducible and its projection over
$V_0\otimes\ldots V_p $ is $D$.
\end{proof}              

\begin{defin}
\label{defdet}
  \rm The hyperdeterminant of $A\in V_0\otimes\ldots\otimes V_p$ is the usual
  determinant of ${\partial}_A$, that is
  \begin{equation}
    \label{eq:partial}
Det(A):=det {\partial}_A    
  \end{equation}
where ${\partial}_A=H^0(\phi_A)$ and $\phi_A:V_0^\ch \otimes
\OO_X\rig{\phi_A}V_1\otimes\OO_X(1,\ldots ,1)$ is the sheaf morphism
associated to $A$.
\end{defin} 
This is  theorem 3.3 of chapter 14 of \cite{GKZ}.
Now, applying remark \ref{rem:det}, we have a
  $GL(V_0)\times \ldots \times GL(V_p)$-equivariant function
\[ Det \colon V_0\otimes\ldots\otimes V_p \to \bigotimes_{i=0}^p\left(det V_i\right)^{\frac{N}{k_i+1}}\]
\[ A\mapsto det ({\partial}_A)\]

\begin{cor}
Let $A\in V_0\otimes\ldots\otimes V_p$ of boundary format. $A$ is nondegenerate
if and only if $Det (A)\neq 0$
\end{cor}
\begin{remark}
  Equality (\ref{eq:partial}) is now proved without any ambiguity of
  the sign, while other methods give an answer modulo 
  the choice of the sign (see \cite{WZ} remark 7.2a). 
\end{remark}

\begin{remark}
\label{perm}
  \rm Any permutation of the $p$ numbers $k_1, \ldots, k_p$ gives
different $m_i$'s and hence a
  different map ${\partial}_A$. As noticed by Gelfand, Kapranov and 
Zelevinsky, in all cases the determinant of ${\partial}_A$ is the
  same by  theorem \ref{main}. 
\end{remark}  
\begin{remark}
\label{notone}
  \rm The given definition of hyperdeterminant can be generalized to
  other cases where the codimension of the degenerate matrices is
bigger than one, these cases are not covered in \cite{GKZ}. 
If $k_0, \ldots, k_p$ are
  nonnegative integers satisfying $k_0=\sum_{i=1}^pk_i$ then we denote
  again $m_j=\sum_{i=1}^{j-1}k_i$ with the convention $m_1=0$.  
  
  Assume we have vector spaces $V_0,\ldots ,V_p$ and a positive
  integer $q$ such that $\dim V_0=q(k_0+1), \dim V_1=q(k_1+1)$ and $\dim
  V_i=(k_i+1)$ for $i=2,\ldots ,p$.  Then the vector spaces
  $V_0\otimes S^{m_1}V_1\otimes\ldots\otimes S^{m_p}V_p$ and $
  S^{m_1+1}V_1\otimes\ldots\otimes S^{m_p+1}V_p$ still have the same
  dimension. In this case degenerate matrices form a subvariety of
  codimension bigger than $1$.
  
  The case $q=p=2$ has been explored in \cite{CO} leading to the proof
  that the moduli space of instanton bundles on $\PP^3$
is affine.
\end{remark}
\section{The Binet-Cauchy theorem for hyperdeterminants of boundary format}
Let $A=(a_{i_0, \dots, i_p})$ a matrix of format $(k_0+1)\times \dots
\times (k_p +1)$ and $B=(b_{j_0, \dots, j_q})$ of format
$(l_0+1)\times \dots \times (l_q +1)$, if $k_p=l_0$ it is defined (see
\cite{GKZ}) the
convolution (or product) $A \ast B$ of $A$ and $B$ as the
$(p+q-1)$-dimensional matrix $C$ of format $(k_0+1)\times \dots \times
(k_{p-1} +1)\times(l_1+1)\times \dots \times (l_q +1)$ with entries
\begin{equation*}
  c_{i_{0},\dots,i_{p-1},j_1,\dots,j_q}=\sum_{h=0}^{k_p}a_{i_0, 
\dots,i_{p-1},h}b_{h,j_1,\dots, j_q}.
\end{equation*}
Similarly, we can define the convolution $A {\ast}_{r,s} B$ with
respect to a pair of indices $r,s$ such that $k_r=l_s$.
\begin{prop}\cite{GKZ}
  If $A,B$ are degenerate then $A \ast B$ is also degenerate and if
  the hyperdeterminants of $A$, $B$ and $A\ast B$ are non-trivial there
  exist polynomials $P(A,B)$ and $Q(A,B)$ in entries of $A$ and $B$
  such that
  \begin{equation*}
    Det(A \ast B)=P(A,B)Det(A)+Q(A,B)Det(B)
  \end{equation*}
\end{prop}
In what follows we prove that in the case of boundary format
matrices the hyperdeterminant of the convolution can be explicitly
described  by only the hyperdeterminants of the involved matrices.
\begin{thm}
  \label{binet}
  If $A\in \Vop$ and $B\in \Woq$ are nondegenerate boundary format
  matrices with $dim V_i=k_i+1$, $dim W_j=l_j+1$ and $W_0^\ch \simeq
  V_p$ then $A\ast B$ is also nondegenerate and
\begin{equation}
  \label{eq:binet}
Det(A\ast
B)=(DetA)^{\binom{l_0}{l_1,\dots,l_q}}(DetB)^{\binom{k_0+1}{k_1,\dots,k_{p-1},k_p+1}}
\end{equation}
\end{thm}
We remark that equation (\ref{eq:binet}) generalizes the Binet-Cauchy 
theorem for determinant of usual square matrices.    
\begin{proof}
  We first observe that the convolution of boundary format matrices
  $A$ and $B$ is also boundary format, then by theorem
  \ref{main} its hyperdeterminant is the usual determinant of
  $\partial_{A\ast B}$

We put 
 \[X_1:=\PVop ;\quad  X_2:=\PWoq \quad \text{and} \quad X:=X_1 \times X_2\]

Since $A$ and $B$ are nondegenerate matrices, they define vector bundles $S_A$ and $S_B$ respectively over $X_1$ and $X_2$ which verify  the following exact sequences
\begin{equation*}
  0\rig{}S_A^\ch \rig{}V_0^\ch\otimes \OO_{X_1}\rig{\phi_A}V_p\otimes
\OO_{X_1}(1,\ldots ,1)\rig{}0. 
\end{equation*}
\begin{equation*}
 0\rig{}S_B^\ch \rig{}W_0^\ch\otimes \OO_{X_2}\rig{\phi_B}W_q\otimes
\OO_{X_2}(1,\ldots ,1)\rig{}0. 
\end{equation*}
Moreover the matrix $A\ast B$ defines the sheaf morphism 
\begin{equation*}
  \begin{CD}
V_0^\ch \otimes \OO_X @>{\phi_{ A\ast B}}>> W_q\otimes\OO_X(1,\ldots ,1)@.    
  \end{CD}
\end{equation*}
If the maps 
\begin{math}
\alpha: X_1\times X_2 \to X_1 \quad \text{and}  \quad \beta: X_1 \times X_2 \to X_2
\end{math}
are the natural projections and $S_{A\ast B} =Ker(\phi_{A\ast B})$, we can construct the following commutative diagram:

\begin{small}
    \begin{equation}
      \begin{CD}
        {} @.{} @. 0 @. 0 @. {} \\
        @. @. @VVV @VVV @. \\
        0 @>>>\alpha^* S_A^\ch  @>{f}>> S_{A\ast B}^\ch @>{g}>> \beta^* S_B^\ch(\underbrace{1,\dots,1}_{p-1},0,\dots0) @. {} \\
        @. @| @VVV @VVV @. \\
        0 @>>> 
        \alpha^* S_A^\ch @>>>
        V_0^\ch \otimes \OO_{X} @>{\alpha^*\phi_A}>>
        V_p \otimes \OO_X(\underbrace{1,\dots,1}_{p-1},0,\dots0) @>>> 0 \\
        @. @. @VV{\phi_{A\ast B}}V @VV{\beta^*\phi_B \otimes
        id_{\OO_X(1,\dots,1,0,\dots,0)}}V @. \\
        {} @.   @.
        W_q \otimes \OO_X(1,\dots, 1) @= W_q \otimes \OO_X(1,\dots, 1)  \\
        @. @. @. @VVV @. \\
        {} @. @. @. 0 @. {}
      \end{CD}
    \end{equation}
  \end{small}
The surjectivity of maps $\beta^*\phi_B \otimes id_{\OO_X(1,\dots,1,0,\dots,0)} $ and $\alpha^*\phi_A$ induce
the surjectivity of $g$ and $\phi_{A\ast B}$, thus $A\ast B$ is
nondegenerate and  $S_{A\ast B}$ is a vector bundle 
   
Moreover, since \[\phi_{A\ast B}= \beta^* \phi_B \otimes
id_{\OO_X(1,\dots,1,0,\dots,0)} \circ \alpha^* \phi_{A}\]
\[\text{and} \qquad \partial_{A\ast B}= H^0 (\phi_{A\ast B} \otimes
id_{\OO(\sum_2^p k_i,\sum_3^p k_i,\dots,k_p,\sum_2^q l_j, \sum_3^q
  l_j, \dots, l_q)})\] then \[\partial_{A\ast B}=(\partial_B \otimes
id_{S^{\sum_2^p k_i+1} V_1 \otimes \dots \otimes S^{k_p+1}V_{p}}
)\circ(\partial_A \otimes id_{S^{\sum_2^q l_j} W_1 \otimes \dots
  \otimes S^{l_q}W_{q-1}})\] i.e. by remark \ref{rem:det}
\[det(\partial_{A\ast B})=(det(\partial_{A}))^{\binom{l_0}{l_1,\dots,l_q}}(det(\partial_{B}))^
{\binom{k_0+1}{k_1,\dots,k_p+1}}\]
as we wanted.
\end{proof}
\begin{remark}
  The degree
  of the hyperdeterminant of a boundary format $(k_0+1)\times \dots
  \times (k_p+1)$ matrix $A$  is given by
  the multinomial coefficient: 
  \begin{equation*}
    \label{eq:deg}
  N_A={\binom{k_0+1} 
  {k_1,\dots,k_p}}
\end{equation*}
This follows also from (\ref{defdet}).
Thus, (\ref{eq:binet}) can be rewritten as 
\begin{equation*}
  Det(A\ast B)={\lbrack (DetA)^{N_B} (DetB)^{N_A} \rbrack}^{\frac{1}{l_0+1}}
\end{equation*}
\end{remark}
\begin{remark}
  The same result of the above theorem works for the convolution
  with respect to the pair of indices $(j,0)$ with $j$ varying in
  $\{1,\dots p \}$. Indeed the condition $W_0^\ch \simeq V_j$ ensure
  that   $A*_{j,0}$ is again of boundary format and  we can arrange the
  indices  as in the proof  because for any permutation
  $\sigma$ we have  $Det(A)=Det(\sigma A)$.
\end{remark}
\begin{cor}
\label{nondeg}
  If $A$ and $B$ are boundary format matrices then 
\[A \ \text{and} \  B \ \text{are nondegenerate} \iff A\ast_{0,j}
B \ \text{are nondegenerate} \]
\end{cor}
The implication $\Longleftarrow$ of the previous corollary is true 
without the assumption of boundary format, see  proposition 1.9 of
\cite{GKZ}.
\begin{remark}
\label{rem:contro}
 Theorem \ref{binet} and the implication $\Longrightarrow$ of
 the corollary \ref{nondeg}
work only for boundary format matrices. Indeed, if, for instance, $A$ and $B$ are $2 \times 2 \times 2$ matrices with 
 \begin{eqnarray*}
   a_{ijk}=0& \ \text{for all} \quad (i,j,k) \notin \{(0,0,0),(1,1,1) \}&  \text{and} \\
   b_{krs}=0& \ \text{for all} \quad (k,r,s) \notin \{(0,0,1),(1,1,0) \}&
\end{eqnarray*}

then $A$ and $B$ are nondegenerate since, applying Cayley formula 
(see \cite{Cay} pag.89 or \cite{GKZ} pag.448), their hyperdeterminants are respectively:
\[
  Det(A)=a_{000}^2a_{111}^2 \quad \text{and} \quad  Det(B)=b_{001}^2b_{110}^2 
\]
but the convolution $A\ast B$ is degenerate. In this case, by using 
Schl\"afli's method of computing hyperdeterminant (\cite{GKZ}), 
it easy to find that $Det(A \ast B)$ corresponds to the 
discriminant of the polynomial 
$F(x_0,x_1)=a_{000}^2a_{111}^2b_{001}^2b_{110}^2x_0^2x_1^2$
which obviously vanishes.  
\end{remark}
\label{biblio}

\vspace{1cm}
\end{document}